\documentclass[12pt]{article}
\usepackage[utf8]{inputenc}
\usepackage{amsthm,amsmath,amsfonts, amssymb}
\usepackage[left=2.25cm,right=2.25cm,top=2.5cm,bottom=2.5cm,
            textwidth=165mm,textheight=250mm,
            nohead, nofoot, footskip=7mm]{geometry}

\usepackage{hyperref}    

\usepackage{authblk}

\hypersetup{
    colorlinks=true,
    linkcolor=blue,
    citecolor=blue,
    filecolor=magenta,      
    urlcolor=cyan,
    pdftitle={The Neumann problem for the generalized H\'enon equation},
    pdfpagemode=FullScreen,
}
\newtheorem{Predl}{\noindent{\bf\textsc{Proposition}}}
\newtheorem{Lemma}{\noindent{\bf\textsc{Lemma}}}
\newtheorem{Theorem}{\noindent{\bf\textsc{Theorem}}}

\newcommand{\Dok}{\par\noindent\textsc{Proof. }}

\newcounter{zam}[section]
\newcommand{\Zam}[1]{\refstepcounter{zam}%
\par\textsc{Remark \arabic{section}.\arabic{zam}. }#1}

\newcommand{\Qa}{{\cal Q}_{p,q,\alpha}}

\newcommand{\minv}{{\rm\bf v}}
\newcommand{\valpha}{{\rm\bf v}_\alpha}
\newcommand{\tpn}{t_{p,n}}
\newcommand{\epn}{\tau_{p,n}}

\title{The Neumann problem for the generalized\\ H\'enon equation. Local analysis}

\author[1]{Alexander I. Nazarov\footnote{Supported by the Ministry of Science and Higher Education of the Russian Federation (agreement 075-15-2025-344 dated 29/04/2025 for Saint Petersburg Leonhard Euler International Mathematical Institute)}}
\author[2]{Alexandra P. Shcheglova}
\affil[1]{St. Petersburg Dept. of Steklov Institute\protect \\ \textit{al.il.nazarov@gmail.com}}
\affil[2]{St. Petersburg Electrotechnical University ``LETI''\protect \\ \textit{alexandra.scheglova@gmail.com}}

\date{}

\begin{document}

\maketitle

\begin{quote}
\noindent\small\textbf{Abstract.} For the boundary value problem 
$$\left\{
\begin{array}{rcll}
-\Delta_p u+u^{p-1}&=&|x|^{\alpha}u^{q-1}&\mbox{in }\Omega,\\
\frac{\displaystyle\partial u}{\displaystyle\partial{\bf n}}&=&0&\mbox{on }\partial \Omega,
\end{array}\right.
$$
in the unit ball $\Omega$, we investigate the properties of the positive radial solution. It is known,  that for $1<p<n$, $\frac{(n-1)p}{n-p}<q<\frac{np}{n-p}$ and sufficiently large $\alpha$ this solution does not provide global minimum to the corresponding energy functional, see~\cite{GazSer} for $p=2$ and~\cite{Shch18} in general case. Nevertheless, it is shown in~\cite{GazSer} that for $n\ge 4$, $p=2$, $2<q<\frac{2n}{n-2}$ and sufficiently large $\alpha$ the radial solution is at least a local minimizer of the energy functional. 

We partially generalize this result. Namely, let $n\ge4$ and let $p>2$ be sufficiently close to $2$. Then for all $p<q<\frac{np}{n-p}$, for sufficiently large $\alpha$ the second variation of the energy functional is positive. The same holds true for all $2<p<n$ if $q>p$ is sufficiently close to $p$. 

\end{quote}
\medskip

\section{Introduction and statement of the problem}

Let $\Omega$ be the unit ball in $\mathbb R^n$ ($n\ge 3$), $S=\partial\Omega$, $p>1$, $q>1$ and $\alpha>0$. We consider the Neumann problem for the generalized H\'enon equation
\begin{equation}
\label{eq:01}
\left\{
\begin{array}{rcll}
-\Delta_p u+u^{p-1}&=&|x|^{\alpha}u^{q-1}&\mbox{in }\Omega,\\[4pt]
u&>&0&\mbox{in }\Omega,\\[4pt]
\frac{\displaystyle\partial u}{\displaystyle\partial{\bf n}}&=&0&\mbox{on }S,
\end{array}
\right.
\end{equation}
where $\Delta_p u={\rm div}(|\nabla u|^{p-2}\nabla u)$ and ${\bf n}$ is the outward unit normal to $S$. 

The Dirichlet problem
\begin{equation}
\label{eq:Henon}
\left\{
\begin{array}{rcll}
-\Delta_p u&=&|x|^\alpha u^{q-1}&\mbox{in}\ \Omega,\\ 
u&>&0,\quad u\Bigl|_S=0&
\end{array}
\right.
\end{equation}
was introduced in 1973 by M. H\'enon~\cite{Hen} for $q>p=2$ as a model for spherically symmetric stellar clusters. In~\cite{Hen} it was studied numerically for some values of $q$ and $\alpha$.

Later the problem~\eqref{eq:Henon} has been extensively studied. Radial solutions were investigated in~\cite{Ni,SSW} for $p=2$ and in~\cite{Naz} for arbitrary $p>1$. The multiplicity of positive solutions under certain conditions on $q$ and $\alpha$ was established in~\cite{SSW} for $p=2$ and in~\cite{NazKol} for all $p>1$.

For the Neumann problem~\eqref{eq:01} with $p=2$, uniqueness (hence radiality) of the positive solution was proved in~\cite{GazSer} for $q>2$ close to $2$ and sufficiently large $\alpha$. Furthermore, the existence of non-radial solutions was established for $\frac{2(n-1)}{n-2}<q<\frac{2n}{n-2}$ and large $\alpha$.  

In the paper~\cite{Shch18}, the existence of non-radial solutions was established for arbitrary $1<p<n$, $\frac{(n-1)p}{n-p}<q<\frac{np}{n-p}$ and large $\alpha$. Moreover, it was proved that for $n\ge 4$ there exist any prescribed number of non-equivalent (not resulting from each other by rotation) positive solutions if $\alpha$ is sufficiently large. Solutions of the problem~\eqref{eq:01} are understood in the weak sense.

It was also shown in~\cite{GazSer} that the radial solution provides a local minimum to the energy functional
\begin{equation}
\Qa(u)=\dfrac{\displaystyle\int\limits_\Omega\bigl(|\nabla u|^p+|u|^p\bigr)dx}
{\left(\displaystyle\int\limits_\Omega |x|^{\alpha}|u|^q dx\right)^{p/q}},\quad u\in W^1_p(\Omega)\setminus\{0\},
\label{eq:02}
\end{equation}
for $n\ge 4$, $p=2$,  $2<q<\frac{2n}{n-2}$ and sufficiently large $\alpha$. Notice that for $q>\frac{2(n-1)}{n-2}$ this minimum is not global.

In the present paper, we establish similar results for $2< p< n$. We stress that in the case $p=2$, the behavior of the radial minimizer is described in terms of Bessel functions. However, for $p\ne 2$ this technique does not work, and the analysis becomes significantly more complicated.

In what follows, the subscript of $o(1)$ indicates the variable with respect to which the smallness is understood. All constants depending only on $n$, $p$ and $q$ are denoted by $C$.

The paper is organized as follows. 
Section 2 is devoted to the generalized Steklov problem, which arises as the limiting problem when $\alpha\to\infty$. In Section 3, we present the main results on the positivity of the second variation of the energy functional, Theorem~\ref{th:main}.\footnote{Unfortunately, we cannot prove a stronger fact, positive definiteness of the second variation. So, the local minimality of the radial solution of~\eqref{eq:01} is still an open problem.} Auxiliary technical Lemmata are collected in the Appendix.

We are grateful to V.E. Bobkov for the useful advice. We also thank anonymous referees (human and artificial) for valuable comments which helped us to improve the exposition.

\section{The Steklov type problem for $p$-Laplacian}

In this section we consider the following problem:
\begin{equation}
\label{eq:Steklov}
\left\{
\begin{array}{rcll}
-\Delta_p u+|u|^{p-2}u&=&0&\mbox{in }\Omega,\\[4pt]
|\nabla u|^{p-2}\,\dfrac{\partial u}{\partial \bf n} &=&\lambda|u|^{p-2}u&\mbox{on }S.
\end{array}
\right.
\end{equation}
It is well known (see, e.g.,~\cite[Theorem 1.1]{MR}) that for any $p>1$ the first eigenvalue $\lambda_p$ of this problem is simple and equals
\begin{equation}
\label{eq:lambda_p}
\lambda_p=\min\limits_{u\in W^1_p(\Omega)\atop u\ne 0}
Q_p(u),\qquad Q_p(u)=\frac{\|u\|^p_{W^1_p(\Omega)}}{\|u\|^p_{L_p(S)}}.
\end{equation}
The corresponding eigenfunction $\varphi_p$ is positive and radial in $\Omega$. 

Denote by $(r,\theta)$ the spherical coordinates of a point $x$. Then $\varphi_p(r)$ is a solution of the problem 
\begin{equation}
\label{eq:varphi_eq}
\left\{
\begin{array}{rr}
-\bigl(r^{n-1}(\varphi_p ')^{p-1}\bigr)'+r^{n-1}\varphi_p^{p-1}=0,&\quad r\in(0,1),\vphantom{\Bigr)} \\ 
\varphi_p '(1)=\lambda_p^{\frac 1{p-1}}\varphi_p(1),& 
\end{array}
\right.
\end{equation}
which is bounded at the origin. We normalize $\varphi_p$  by the assumption $\|\varphi_p\|_{W^1_p(\Omega)}=1$, which is equivalent to
\begin{equation*}
\varphi_p(1)=\lambda_p^{-\frac{1}{p}}({\rm meas}\,S)^{-\frac{1}{p}}.
\end{equation*}
Standard argument shows that $\varphi_p(r)$ is strictly increasing and has the following asymptotic behavior at zero:
\begin{equation}\label{eq:varphi}
\varphi_p(r)=\varphi_0+o_r(1),\qquad \varphi_p'(r)=\varphi_0n^{-\frac{1}{p-1}}r^{\frac{1}{p-1}}\bigl(1+o_r(1)\bigr),\quad r\to +0.
\end{equation}
Here the positive constant $\varphi_0=\varphi_0(n,p)$ is determined by the normalization conditions. \medskip

Next, for $p\ge 2$ we consider the second variation of the functional~\eqref{eq:lambda_p} at the point $\varphi_p$. We decompose an arbitrary increment $h\in W^1_p(\Omega)$ as
\begin{equation}
\label{eq:h}
h(r,\theta)=h_1(r)+g(r,\theta),\quad\text{where }h_1(r)=({\rm meas}\,S)^{-1}\int\limits_S h(r,\theta)dS,\quad g\in{\cal W},
\end{equation}
$$
{\cal W}=\left\{g\in W^1_2(\Omega):\int\limits_S g(r,\theta)dS=0\text{ for a.e. }r\in(0,1)\right\}.
$$

A direct calculation yields
$$
D^2 Q_p(\varphi_p,h)=D^2 Q_p(\varphi_p,h_1)+p\lambda_p\Bigl(F_p (g)-\lambda_p(p-1)\varphi_p^{p-2}(1)\Bigr)\cdot\int\limits_S g^2(1,\theta)dS,
$$
where
\begin{equation}
F_p (g)=\frac{\displaystyle\int\limits_{\Omega}\Bigl(
(p-2)|\nabla\varphi_p|^{p-4}\langle\nabla\varphi_p,\nabla g\rangle^2+
|\nabla\varphi_p|^{p-2}|\nabla g|^2+
(p-1)\varphi_p^{p-2} g^2\Bigr)dx}
{\displaystyle\int\limits_{S}g^2(1,\theta)dS}.
\label{eq:Fpg}
\end{equation}

\begin{Lemma}
\label{lm:2}
We have
\begin{equation}
\label{eq:min_F_p}
\min_{g\in{\cal W}, g(1,\theta)\not\equiv 0}F_p(g)=
\lambda_p^{\frac 2p-\frac 1{p-1}-1}({\rm meas}\,S)^{\frac 2p -1}
\bigl(1-(n-1)\lambda_p\bigr).
\end{equation}
\end{Lemma}

\Dok{The attainability of the minimum in~\eqref{eq:min_F_p} easily follows by compactness argument.  The minimization problem admits separation of variables, so we can write $g(r,\theta)=w(r)\psi(\theta)$ with $\displaystyle\int\limits_S\psi(\theta)dS=0$. Then
$$
F_p(g)=(p-1)w^{-2}(1)\int\limits_0^1\bigl((\varphi_p')^{p-2}w'^2+\varphi_p^{p-2}w^2\bigr)r^{n-1}dr
+w^{-2}(1)\int\limits_0^1(\varphi_p')^{p-2}w^2 r^{n-3}dr\cdot\dfrac{\int\limits_S|\nabla_\theta\psi|^2dS}{\int\limits_S\psi^2dS},
$$
where $\nabla_\theta$ is the tangential gradient on $S$.

The minimum of $\dfrac{\int\limits_S|\nabla_\theta\psi|^2dS}{\int\limits_S\psi^2dS}$ over functions with zero mean equals $(n-1)$, attained at the first-order spherical harmonics. Thus,
\begin{multline*}
\min F_p(g)=\min_{w\in W^1_p\left((0,1),r^{n-1}dr\right),\atop w(1)\ne0}\frac1{w^2(1)}\Biggl[ (p-1)\int\limits_0^1\bigl(|\varphi_p'|^{p-2}w'^2 
+\varphi_p^{p-2}w^2\bigr)r^{n-1}dr
\\+(n-1)\int\limits_0^1|\varphi_p'|^{p-2}w^2 r^{n-3}dr\Biggr].
\end{multline*}

The Euler-Lagrange equation for this problem is
\begin{equation}
\label{eq:varphi_prime}
\left\{
\begin{array}{ll}
(p-1)\Bigl(\bigl(-r^{n-1}(\varphi_p')^{p-2}w'\bigr)'+r^{n-1}\varphi_p^{p-2}w\Bigr)+(n-1)r^{n-3}(\varphi_p')^{p-2}w=0,& r\in(0,1),\\[4pt]
(p-1)(\varphi_p'(1))^{p-2}w'(1)=\lambda w(1),
\end{array}\right.
\end{equation}
where $\lambda$ is the Lagrange multiplier. By differentiating~\eqref{eq:varphi_eq} we verify that $w=\varphi_p'$ satisfies the first equation in~\eqref{eq:varphi_prime}. Substituting $r=1$ into~\eqref{eq:varphi_eq} and using the boundary condition we obtain
$$
(p-1)w'(1)=\left(\frac1{\lambda_p}-(n-1)\right)w(1).
$$
Hence
\begin{multline*}
\lambda= \min F_p(g)=(\varphi_p'(1))^{p-2}\left(\frac1{\lambda_p}-(n-1)\right)
=\lambda_p^{\frac{p-2}{p-1}-1}\varphi_p(1)^{p-2}(1-(n-1)\lambda_p)\\
=\lambda_p^{\frac2p-\frac1{p-1}-1}({\rm meas}\,S)^{\frac2p-1}(1-(n-1)\lambda_p).
\end{multline*}

Finally, it is easy to see that any solution of equation~\eqref{eq:varphi_prime} linearly independent with $\varphi'_p(r)$ has asymptotics $\widetilde\varphi(r)\asymp r^{1-n}$ as $r\to +0$.  Since $\widetilde\varphi\not\in W^1_p(\Omega)$, we conclude $w=C\varphi'_p$. \hfill$\square$}

\section{The properties of radial solutions of problem~\eqref{eq:01}}

Repeating verbatim the proof of Proposition~1.1 in~\cite{Naz}, we obtain:

\begin{Predl}
For $q\le p$ problem~\eqref{eq:01} has a unique solution, and it is radial.
\end{Predl}

Therefore, in what follows we assume $q>p$.\medskip

Solutions of problem~\eqref{eq:01} can be obtained as critical points of the energy functional~\eqref{eq:02}.
Indeed, the Euler--Lagrange equation for $\Qa$ transforms into~\eqref{eq:01} by a suitable rescaling, provided that $u>0$.\smallskip

For $1<p<n$ we denote by $p^*=\frac{np}{n-p}$ the critical Sobolev exponent. The functional $\Qa$ is well defined for all $q\in(p,p^*)$. Its restriction to the subspace $W^1_{p,rad}(\Omega)$ of radially symmetric functions is well-defined for all $q\in\left(p,p^*+\frac{p\alpha}{n-p}\right)$. 

\begin{Predl}[{\cite[Prop. 2.1]{Shch18}}] 
Let $\alpha>0$ and $q\in\left(p,p^*+\frac{p\alpha}{n-p}\right)$. Then $\Qa$ attains a positive minimum on $W^1_{p,rad}(\Omega)$, and the minimizer $\minv_\alpha$ (after multiplication by a suitable constant) is a positive solution of problem~\eqref{eq:01}.
\end{Predl}

Without loss of generality we assume $\valpha$ normalized by $\|\minv_\alpha\|_{W^1_p(\Omega)}=1$. 

For each $\alpha>0$ and $q\in\left(p,p^*+\frac{\alpha p}{n-p}\right)$, denote
\begin{equation}
\mu_{p,q,\alpha}=\min_{v\in W^1_{p,rad}(\Omega), v\ne 0}\Qa(v)=\left(\mathrm{meas}\, S\int\limits_0^1 r^{\alpha+n-1}\valpha^q dr\right)^{-p/q}.
\label{eq:03_0}
\end{equation}
Then $\minv_\alpha$ is a radial solution of the problem
\begin{equation}
\label{eq:03}
\left\{
\begin{array}{rcll}
-\Delta_p u+u^{p-1}&=&\mu_{p,q,\alpha}^{q/p}|x|^{\alpha}u^{q-1}&\mbox{in }\Omega,\\[4pt]
\frac{\displaystyle\partial u}{\displaystyle\partial{\bf n}}&=&0&\mbox{on }\partial\Omega
\end{array}
\right.
\end{equation}
which can be rewritten as follows:
\begin{equation}
\label{eq:03_1}
\left\{
\begin{array}{rcll}
-\bigl(r^{n-1}|\valpha'|^{p-2}\valpha'\bigr)'+r^{n-1}\valpha^{p-1}&=&\mu_{p,q,\alpha}^{q/p}r^{\alpha+n-1}\valpha^{q-1},& r\in(0,1),\\[4pt]
\valpha'(1)&=&0.&
\end{array}
\right.
\end{equation}

Problem~\eqref{eq:Steklov} serves as the limit problem for~\eqref{eq:03} as $\alpha\to+\infty$ in the following sense.

\begin{Theorem}[{\cite[Theorems 2.1, 2.2]{Shch18}}]
\label{th:lim}
Let $p>1$ and $q\in(p,+\infty)$. Then, as $\alpha\to+\infty$:
\begin{itemize}
\item[$(a)$] $\valpha\to\varphi_p$ in $W^1_p(\Omega)$;
\item[$(b)$] $\valpha\to\varphi_p$ in ${\cal C}(\overline\Omega)$;
\item[$(c)$] $\valpha$ are uniformly bounded in ${\cal C}^1(\overline\Omega)$, and
\begin{equation}
\label{eq:v_alpha_prime}
\valpha'(r)=\varphi_p'(r)\bigl(1+o_r(1)\bigr)\bigl(1+o_\alpha(1)\bigr),\quad r\to +0;
\end{equation}
\item[$(d)$] for any $\delta\in(0,1)$, $\valpha\to\varphi_p$ in ${\cal C}^1(\overline\Omega_\delta)$, where $\Omega_\delta$ is the ball of radius $\delta$;
\item[$(e)$] $\mu_{p,q,\alpha}=(\alpha+n)^{p/q}\bigl({\rm meas}\,S\bigr)^{1-p/q}\lambda_p(1+o_\alpha(1))$.
\end{itemize}
\end{Theorem}

In what follows we need also some new estimates.

\begin{Lemma}
\label{lm:v_alpha_1}
Let $p>1$ and $q\in(p,+\infty)$. Then for sufficiently large $\alpha$ we have $\valpha'(r)>0$. Moreover, there exists $M=M(n,p,q)$ such that, as $\alpha\to +\infty$ we have 
\begin{align}
\label{eq:asymp_3}
Cr  \le (\valpha'(r))^{p-1}\le\, & Cr, && r\in(0,1-\tfrac{M}{\alpha}],\\
(\valpha'(r))^{p-1} \ge\, & C\dfrac{\alpha}{\ln\alpha}(1-r), && r\in\left(1-\tfrac{3\ln\alpha}{2\alpha},1-\tfrac{M}{\alpha}\right],
\label{eq:asymp}
\\
C\alpha(1-r)  \le (\valpha'(r))^{p-1}\le\, & C\alpha(1-r), && r\in\left(1-\tfrac{M}{\alpha},1\right).
\label{eq:asymp_1}
\end{align}
\end{Lemma}
\Dok{The right part of~\eqref{eq:asymp_3} follows from~\eqref{eq:varphi} and~\eqref{eq:v_alpha_prime}. To prove the left part, we integrate the equation~\eqref{eq:03_1} from $0$ to $r$. The substitution at~$0$ vanishes due to~\eqref{eq:v_alpha_prime}, and for sufficiently large $\alpha$ we have by Theorem~\ref{th:lim}$(e)$
\begin{multline}
\label{eq:asymp_2}
r^{n-1}|\valpha'(r)|^{p-2}\valpha'(r)
 =\int\limits_0^r \valpha^{p-1}(s) s^{n-1} ds-\mu_{p,q,\alpha}^{q/p}\int\limits_0^r s^{\alpha+n-1}\valpha^{q-1}(s)\,ds \ge
 Cr^n-\dfrac{C\mu_{p,q,\alpha}^{q/p}}{\alpha+n}r^{\alpha+n}\\
 \ge Cr^n\left(1-C\left(1-\tfrac{M}{\alpha}\right)^\alpha\right)\ge Cr^n(1-Ce^{-M})
\end{multline}
if $r\in(0,1-\frac{M}{\alpha}]$. The lattely bracket in~\eqref{eq:asymp_2} is greater than $\frac{1}{2}$ for some $M=M(n,p,q)$ and we get the left inequality in~\eqref{eq:asymp_3}.


Next, let $r\in(1-\frac{3\ln\alpha}{2\alpha},1)$. We integrate the equation~\eqref{eq:03_1} from $r$ to $1$. This gives for sufficiently large $\alpha$ 
\begin{multline*}
r^{n-1}|\valpha'(r)|^{p-2}\valpha'(r)
 =\mu_{p,q,\alpha}^{q/p}\int\limits_r^1 s^{\alpha+n-1}\valpha^{q-1}(s)\,ds-\int\limits_r^1 \valpha^{p-1}(s) s^{n-1} ds \ge
 C(1-r^{\alpha+n})-C(1-r^n)\\
 \ge C(1-r)\left(\frac{1-r^{\alpha+n}}{1-r}-C\right)\ge C(1-r)\left(\dfrac{\alpha}{\ln\alpha}-C\right)>0,
\end{multline*}
and the positivity of $\valpha'(r)$ is proved for all $r\in(0,1)$. Also this provides~\eqref{eq:asymp}.

Finally, for $r\in(1-\frac{M}{\alpha},1)$ we again integrate~\eqref{eq:03_1} from $r$ to $1$ and apply Theorem~\ref{th:lim}$(b,e)$:
\begin{multline*}
r^{n-1}(\valpha'(r))^{p-1} =\int\limits_r^1 \Bigl(C\alpha(\varphi_p(s)+o_\alpha(1))^{q-1}s^{n+\alpha-1}  -C(\varphi_p(s)+o_\alpha(1))^{p-1}s^{n-1}\Bigr)ds \\
=C\bigl(1+o_r(1)\bigr)(1-r^{n+\alpha})-C\bigl(1+o_r(1)\bigr)(1-r^n).
\end{multline*}
By the Lagrange mean value theorem, for some $\hat r\in\left(1-\frac{M}{\alpha},1\right)$
we obtain
$$
1-r^{n+\alpha}=(n+\alpha) \hat r^{\,n+\alpha-1}(1-r)\asymp \alpha(1-r),
$$
whereas $1-r^{n}\asymp 1-r$. This provides~\eqref{eq:asymp_1}. 
\hfill$\square$\medskip
}

The following result from~\cite{Shch18} shows that for certain parameters, the global minimizer cannot be radial.

\begin{Predl}
[{\cite[Theorem 3.1]{Shch18}}]
Let $1<p<n$, $q\in(p^{**},p^*)$. Then there exists $\widehat\alpha(p,q)>0$ such that for all $\alpha>\widehat\alpha$, the global minimizer of $\Qa$ is not radial.
\end{Predl}

Here $p^{**}=\frac{(n-1)p}{n-p}$ is the critical trace Sobolev embedding exponent.

Nevertheless, we show that for appropriate $2<p<n$ and $q>p$, the second variation of ${\cal Q}_{p,q,\alpha}$ at the radial function $\valpha$ is strictly positive for all sufficiently large $\alpha$. For $p=2$ such a result was obtained in~\cite{GazSer}. \medskip

Here we formulate the main results of the paper.

\begin{Theorem}
\label{th:main}
{\bf 1}. Let $2 <p<n$. Then there exists $q_{loc}(n,p)>p$ such that for each $q\in(p,q_{loc})$ and for all sufficiently large $\alpha$, we have $D^2{\cal Q}_{p,q,\alpha}(\valpha;h,h)>0$ for any $h\in W^1_p(\Omega)$ such that $h\ne C\valpha$.

{\bf 2}. Let $n\ge4$. Then there exists $p_{loc}(n)\in(2,n)$ such that for any $p\in(2,p_{loc})$ and $q\in(p,p^*)$, for all sufficiently large $\alpha$ we have $D^2{\cal Q}_{p,q,\alpha}(\valpha;h,h)>0$ for any $h\in W^1_p(\Omega)$ such that $h\ne C\valpha$.
\end{Theorem}

\Dok
By the principle of symmetric criticality (see, e.g.,~\cite{Pal} and~\cite{KapLad}), the first differential $D\Qa(\minv_\alpha;h)$ vanishes for all $h\in W^1_p(\Omega)$, not only for radial increments. We decompose an arbitrary increment $h\in W^1_p(\Omega)$ as it is shown in~\eqref{eq:h}. Direct computation gives
$$
D^2\Qa(\valpha;h,h)=D^2\Qa(\valpha;h_1,h_1)+D^2\Qa(\valpha;g,g).
$$
Since $\valpha$ minimizes $\Qa$ among radial functions, we have $D^2\Qa(\valpha;h_1,h_1)\ge0$. Moreover, it is easy to see that $D^2\Qa(\valpha;\valpha,\valpha)=0$. Thus it suffices to show that under assumptions of the Theorem we have
\begin{eqnarray}
D^2\Qa(\valpha;h_1,h_1)>0 &&\mbox{for any}\  h_1\ne C\valpha,
\label{eq:D1}\\   D^2\Qa(\valpha;g,g)>0 &&\mbox{for any non-zero}\  g\in{\cal W}.
\label{eq:D2}
\end{eqnarray}
These inequalities are proved, respectively, in Subsections~\ref{proof_1} and~\ref{proof_2}. 
\hfill$\square$

\subsection{Proof of inequality~\eqref{eq:D1}}
\label{proof_1}

We take into account that $D\mathcal{Q}_{p,q,\alpha}(\valpha;h_1)=0$. Then direct calculation gives 
\begin{multline*}
\dfrac{1}{p}D^2\Qa(\valpha;h_1,h_1)={\left(\int\limits_\Omega |x|^{\alpha}|\valpha|^q dx\right)^{-p/q}}\times\\
\times\int\limits_\Omega \Bigl((p-2)|\nabla \valpha|^{p-4}\langle\nabla\valpha,\nabla h_1\rangle^2+|\nabla\valpha|^{p-2}|\nabla h_1|^2+(p-1)\valpha^{p-2}h_1^2\Bigr)\,dx\\
+(q-p){\left(\int\limits_\Omega |x|^{\alpha}|\valpha|^q dx\right)^{-p/q-2}}\left(\int\limits_\Omega |x|^\alpha\valpha^{q-1}h_1 dx\right)^2\int\limits_\Omega \bigl(|\nabla\valpha|^p+\valpha^p\bigr)dx\\
-(q-1){\left(\int\limits_\Omega |x|^{\alpha}|\valpha|^q dx\right)^{-p/q-1}}\int\limits_\Omega |x|^\alpha\valpha^{q-2}h_1^2 dx\cdot\int\limits_\Omega \bigl(|\nabla\valpha|^p+\valpha^p\bigr)dx,
\end{multline*}
or, using~\eqref{eq:03_0},
\begin{multline*}
\dfrac{1}{p}D^2\Qa(\valpha;h_1,h_1)=\mu_{p,q,\alpha}(\mathrm{meas}\,S)\left[(p-1)\int\limits_0^1 r^{n-1}\bigl((\valpha')^{p-2}h_1'\,^2+\valpha^{p-2}h_1^2\bigr)dr \right.\\
\left.+(q-p)\mu_{p,q,\alpha}^{2q/p}(\mathrm{meas}\,S)\left(\int\limits_0^1 r^{\alpha+n-1}\valpha^{q-1}h_1 dr\right)^2
-(q-1)\mu_{p,q,\alpha}^{q/p}\int\limits_0^1 r^{\alpha+n-1}\valpha^{q-2}h_1^2 dr
\right].    
\end{multline*}

We proceed by contradiction. Let there exist the function $h_1(r)=w(r)\valpha(r)$, $w(r)\ne const$ such that
\begin{multline*}
0=\dfrac{1}{p\,\mu_{p,q,\alpha}(\mathrm{meas}\,S)} D^2\Qa(\valpha;h_1,h_1)\\
=(p-1)\int\limits_0^1 r^{n-1}\Bigl((\valpha')^{p}w^2+2(\valpha')^{p-1}\valpha w'w+(\valpha')^{p-2}\valpha^2 w'\,^2+\valpha^p w^2\Bigr)dr\\    
+(q-p)\mu_{p,q,\alpha}^{2q/p}(\mathrm{meas}\, S)\left(\int\limits_0^1 r^{\alpha+n-1}\valpha^q w\,dr\right)^2-(q-1)\mu_{p,q,\alpha}^{q/p}\int\limits_0^1 r^{\alpha+n-1}\valpha^q w^2 dr.
\end{multline*}

After multiplying~\eqref{eq:03_1} by $\valpha w^2$ and integrating from $0$ to $1$ we obtain
\begin{equation*}
\int\limits_0^1 r^{n-1}\Bigl((\valpha')^{p}w^2+2(\valpha')^{p-1}\valpha w'w+\valpha^p w^2\Bigr) dr=  \mu_{p,q,\alpha}^{q/p}\int\limits_0^1 r^{\alpha+n-1}\valpha^q w^2 dr  .
\end{equation*}
Hence
\begin{multline*}
0=(p-1)\int\limits_0^1 r^{n-1}(\valpha')^{p-2}\valpha^2 w'\,^2 dr - (q-p)\mu_{p,q,\alpha}^{q/p}\int\limits_0^1 r^{\alpha+n-1}\valpha^q w^2 dr\\
+(q-p)\mu_{p,q,\alpha}^{2q/p}(\mathrm{meas}\, S)\left(\int\limits_0^1 r^{\alpha+n-1}\valpha^q w\,dr\right)^2.
\end{multline*}

Next, we decompose $w(r)=\mathbf{w}(r)+W_0$, where $\mathbf{w}\ne 0$ and
\begin{equation}
W_0=\dfrac{\int\limits_0^1 r^{\alpha+n-1}\valpha^q w\,dr}{\int\limits_0^1 r^{\alpha+n-1}\valpha^q\, dr},\qquad
\int\limits_0^1 r^{\alpha+n-1}\valpha^q \mathbf{w}\,dr=0.
\label{eq:ortogon}
\end{equation}

Finally, we have
\begin{equation*}
0=(p-1)\int\limits_0^1 r^{n-1} (\valpha')^{p-2}\valpha^2\mathbf{w}'\,^2 dr- (q-p)\mu_{p,q,\alpha}^{q/p}\int\limits_0^1 r^{\alpha+n-1}\valpha^q\mathbf{w}^2 dr=:\mathcal{G}_{p,q,\alpha}(\mathbf{w}).
\end{equation*}

Since the functional $\mathcal{G}_{p,q,\alpha}$ is nonnegative, 
$\mathbf{w}$ provides its minimum on the subspace of functions in $W^1_2((0,1),r^{n-1}dr)$ with the orthogonality condition~\eqref{eq:ortogon}. Therefore $\mathbf{w}$ satisfies the Euler--Lagrange equation  
\begin{equation}
-(p-1)\Bigl(r^{n-1}(\valpha')^{p-2}\valpha^2\mathbf{w}'\Bigr)'=(q-p)\mu_{p,q,\alpha}^{q/p}r^{\alpha+n-1}\valpha^q\mathbf{w}+\mu_1 r^{\alpha+n-1}\valpha^q   
\label{eq:w}
\end{equation}
with natural boundary condition
\begin{equation}
\lim\limits_{r\to 1-0} \bigl(\valpha'(r)\bigr)^{p-2}\mathbf{w}'(r)=0.
\label{eq:w_prime_1}
\end{equation}

Due to~\eqref{eq:varphi} and~\eqref{eq:v_alpha_prime} we have, as $r\to+0$,
\begin{equation}
\Bigl(C_1(1+o_r(1))r^{n-1+\frac{p-2}{p-1}}\mathbf{w}'\Bigr)'=C_2(1+o_r(1))\mu_{p,q,\alpha}^{q/p}r^{\alpha+n-1}\mathbf{w}+C_3(1+o_r(1))\mu_1 r^{\alpha+n-1}.
\label{eq:w1}
\end{equation}
We seek the asymptotics of the solution in the form
$$
\mathbf{w}(r)\sim Cr^{\beta},\quad r\to +0.
$$
The order of the left-hand  side of~\eqref{eq:w1} is $n-3+\frac{p-2}{p-1}+\beta$ while the order of the right-hand side is $\alpha+n-1+\min\{\beta,0\}$. For $\alpha$ sufficiently large, the right-hand side can be considered as a small perturbation of equation
\begin{equation*}
\Bigl(C_1(1+o_r(1))r^{n-1+\frac{p-2}{p-1}}\mathbf{w}'\Bigr)'=0.
\end{equation*}
This equation has solutions with asymptotics $\mathbf{w}(r)\sim Cr^\beta$ if 
$$
\beta=0\quad\mathrm{or}\quad\beta=-\Bigl(n-2+\dfrac{p-2}{p-1}\Bigr).
$$
The second solution does not belong to $W^1_2((0,1),r^{n-1}dr)$. Therefore, $\mathbf{w}(r)\asymp 1$. Integrating~\eqref{eq:w} from $0$ to $r$, we obtain $\mathbf{w}'(r)\asymp r^{\alpha-\frac{p-2}{p-1}}$,\quad $r\to +0$. Due to these estimates, we can integrate~\eqref{eq:w} from 0 to 1. Taking~\eqref{eq:w_prime_1} into account, we obtain $\mu_1=0$. So $\mathbf{w}$ satisfies the homogeneous equation
\begin{equation}
-(p-1)\Bigl(r^{n-1}(\valpha')^{p-2}\valpha^2\mathbf{w}'\Bigr)'=(q-p)\mu_{p,q,\alpha}^{q/p}r^{\alpha+n-1}\valpha^q\mathbf{w}.
\label{eq:w_hom}
\end{equation}

By~\cite[Lemma 1]{Str} we have
$$
|\mathbf{w}(r)|\le \dfrac{C\|\mathbf{w}\|_{W^1_2(\Omega)}}{r^{\frac{n-2}{2}}}.
$$
Using Theorem~\ref{th:lim}$(b,e)$ we can estimate the right-hand side of~\eqref{eq:w_hom} for sufficiently large $\alpha$ as
\begin{equation}
\Bigl|(q-p)\mu_{p,q,\alpha}^{q/p}r^{\alpha+n-1}\valpha^q\mathbf{w}\Bigr|\le C\alpha\, r^{\alpha+\frac{n}{2}}\, \|\mathbf{w}\|_{W^1_2(\Omega)}.
\label{eq:w_ineq}
\end{equation}


For $r\in(0,1-\frac{\ln\alpha}{\alpha})$ we integrate~\eqref{eq:w_hom} from 0 to $r$. Using~\eqref{eq:asymp_3} and~\eqref{eq:w_ineq} we get
$$
r^{n-1+\frac{p-2}{p-1}}|\mathbf{w}'(r)|\le C\alpha\int\limits_0^r s^{\alpha+\frac{n}{2}} dr\cdot\|\mathbf{w}\|_{W^1_2(\Omega)}\le C
r^{\alpha+\frac{n}{2}+1}\cdot\|\mathbf{w}\|_{W^1_2(\Omega)},
$$
so,
\begin{equation}
\label{eq:27}
|\mathbf{w}'(r)|\le C r^{\alpha-\frac{n}{2}+\frac{p}{p-1}}\cdot\|\mathbf{w}\|_{W^1_2(\Omega)} 
\end{equation}
and
$$
r^{n-1}|\mathbf{w}'(r)|^2\le C r^{2\alpha+\frac{p+1}{p-1}}\cdot\|\mathbf{w}\|^2_{W^1_2(\Omega)}\le C\left(1-\dfrac{\ln\alpha}{\alpha}\right)^{2\alpha} \cdot\|\mathbf{w}\|^2_{W^1_2(\Omega)}=o_\alpha(1)\cdot \|\mathbf{w}\|^2_{W^1_2(\Omega)}.
$$

For $r\in[1-\frac{\ln\alpha}{\alpha},1)$ we integrate~\eqref{eq:w_hom} from $r$ to 1. Using~\eqref{eq:asymp},~\eqref{eq:asymp_1} and boundary condition~\eqref{eq:w_prime_1}, we obtain
\begin{multline*}
\left(\dfrac{\alpha}{\ln\alpha}\right)^{\frac{p-2}{p-1}}(1-r)^{\frac{p-2}{p-1}}|\mathbf{w}'(r)|\le C\alpha\int\limits_r^1 s^{\alpha+\frac{n}{2}} dr\cdot\|\mathbf{w}\|_{W^1_2(\Omega)} \\ \le C(1-
r^{\alpha+\frac{n}{2}+1})\cdot\|\mathbf{w}\|_{W^1_2(\Omega)}\le C\alpha(1-r)\cdot\|\mathbf{w}\|_{W^1_2(\Omega)},
\end{multline*}
so,
\begin{equation}
\label{eq:28}
|\mathbf{w}'(r)|\le C \alpha^{\frac{1}{p-1}}(\ln\alpha)^{\frac{p-2}{p-1}}(1-r)^{\frac{1}{p-1}} \|\mathbf{w}\|_{W^1_2(\Omega)}
\end{equation}
and
$$
r^{n-1}|\mathbf{w}'(r)|^2\le C \alpha^{\frac{2}{p-1}}(\ln\alpha)^{\frac{2(p-2)}{p-1}}(1-r)^{\frac{2}{p-1}}\cdot\|\mathbf{w}\|^2_{W^1_2(\Omega)}\le C(\ln\alpha)^{2}\cdot \|\mathbf{w}\|^2_{W^1_2(\Omega)}.
$$

Therefore, for sufficiently large $\alpha$
\begin{equation*}
\int\limits_0^1 r^{n-1}\mathbf{w}'\,^2 dr=\Bigl(\!\int\limits_0^{1-\frac{\ln\alpha}{\alpha}} \!+\! \int\limits_{1-\frac{\ln\alpha}{\alpha}}^1\!\Bigr)\, r^{n-1}\mathbf{w}'\,^2 dr
\le \left(o_\alpha(1)  + \dfrac{C(\ln\alpha)^{3}}{\alpha}\right) \|\mathbf{w}\|^2_{W^1_2(\Omega)}=o_\alpha(1)\cdot\|\mathbf{w}\|^2_{W^1_2(\Omega)},
\end{equation*}
that is 
\begin{equation}
\label{eq:29}
\|\nabla \mathbf{w}\|^2_{L_2(\Omega)}=o_\alpha(1)\cdot\|\mathbf{w}\|^2_{W^1_2(\Omega)}.    
\end{equation}

By the relation~\eqref{eq:ortogon}, where exists $r_0\in(0,1)$ such that $\mathbf{w}\bigl|_{r=r_0}=0$. So, using~\eqref{eq:27} and~\eqref{eq:28} we deduce  
\begin{multline*}
|\mathbf{w}(r)|=\left|\int\limits_{r_0}^r \mathbf{w}'(s)\,ds\right|\le \int_0^1|\mathbf{w}'(s)|ds \\
\le \Biggl( C \int\limits_0^{1-\frac{\ln\alpha}{\alpha}}  s^{\alpha-\frac{n}{2}+\frac{p}{p-1}} ds + C \alpha^{\frac{1}{p-1}}(\ln\alpha)^{\frac{p-2}{p-1}}\int\limits_{1-\frac{\ln\alpha}{\alpha}}^1 (1-s)^{\frac{1}{p-1}} ds\Biggr) \|\mathbf{w}\|_{W^1_2(\Omega)}\\
\le\left(\dfrac{C}{\alpha^2}+\dfrac{C(\ln\alpha)^2}{\alpha}\right)\|\mathbf{w}\|_{W^1_2(\Omega)}=o_\alpha(1)\cdot \|\mathbf{w}\|_{W^1_2(\Omega)},
\end{multline*}
and we arrive at
\begin{equation*}
\|\mathbf{w}\|^2_{L_2(\Omega)}=o_\alpha(1)\cdot\|\mathbf{w}\|^2_{W^1_2(\Omega)}.    
\end{equation*}
The latter relation and~\eqref{eq:29} imply $\mathbf{w}=0$ for sufficiently large $\alpha$. This completes the proof of~\eqref{eq:D1}.

\subsection{Proof of inequality~\eqref{eq:D2}}\label{proof_2}

We denote
\begin{multline*}
{\cal F}_{p,q,\alpha}(g):=\dfrac{1}{p} D^2\Qa(\valpha;g,g)=
(p-2)\int\limits_\Omega |\nabla\valpha|^{p-4}\langle\nabla\valpha,\nabla g\rangle^2 dx
+\int\limits_\Omega |\nabla\valpha|^{p-2}|\nabla g|^2 dx \\
+(p-1)\int\limits_\Omega \valpha^{p-2}g^2 dx
-(q-1)\mu_{p,q,\alpha}^{q/p}\int\limits_\Omega |x|^{\alpha}\valpha^{q-2}g^2 dx.
\end{multline*}
We again proceed by contradiction and assume that ${\cal F}_{p,q,\alpha}(g)$ takes non-positive values for arbitrarily large $\alpha$ and some $g=g_\alpha\ne 0$.\smallskip

The functional ${\cal F}_{p,q,\alpha}$ is homogeneous, so we can consider it on the set 
$$
{\cal W}_1=\left\{g\in{\cal W}:\|g\|_{W^1_2(\Omega)}=1\right\}.
$$

\begin{Lemma}
\label{lm:achiev}
If ${\cal F}_{p,q,\alpha}$ takes a non-positive value at some function in ${\cal W}_1$, then the infimum 
$$
\inf_{{\cal W}_1}{\cal F}_{p,q,\alpha}=-\lambda_{p,q,\alpha}\le0
$$
is attained.
\end{Lemma}

\Dok{Let $\{g_k\}\subset{\cal W}_1$ be a minimizing sequence. Assume $g_k\rightharpoondown \widehat{g}$ in $W^1_2(\Omega)$. 
Set $w_k=g_k-\widehat{g}$. Then 
\begin{multline*}
-\lambda_{p,q,\alpha}+o_k(1)={\cal F}_{p,q,\alpha}(\widehat{g}+w_k) = (p-2)\int\limits_\Omega |\nabla\valpha|^{p-4}\langle\nabla\valpha,\nabla (\widehat{g}+w_k)\rangle^2dx \\
+\int\limits_\Omega |\nabla\valpha|^{p-2}|\nabla (\widehat{g}+w_k)|^2dx+(p-1)\int\limits_\Omega \valpha^{p-2}(\widehat{g}+w_k)^2 dx-(q-1)\mu_{p,q,\alpha}^{q/p}\int\limits_\Omega |x|^{\alpha}\valpha^{q-2}(\widehat{g}+w_k)^2 dx\\
\stackrel{(*)}{=}{\cal F}_{p,q,\alpha}(\widehat{g})+(p-2)\int\limits_\Omega |\nabla\valpha|^{p-4}\langle\nabla\valpha,\nabla w_k\rangle^2dx+\int\limits_\Omega |\nabla\valpha|^{p-2}|\nabla w_k|^2dx+o_k(1).
\end{multline*}
The equality $(*)$ follows from the convergence $w_k\rightharpoondown 0$ in $W^1_2(\Omega)$ and $w_k\to 0$ in $L_2(\Omega)$.

Passage to the limit gives 
$$
-\lambda_{p,q,\alpha}\ge \mathcal{F}_{p,q,\alpha}(\widehat{g}).
$$

If $\lambda_{p,q,\alpha}=0$, then it is attained, since $\mathcal{F}_{p,q,\alpha}$ takes a non-positive value on $\mathcal{W}_1$. Otherwise $\widehat{g}\ne 0$, and we have 
\begin{equation*}
-\lambda_{p,q,\alpha}\ge \|\widehat{g}\|^2_{W^1_2(\Omega)}{\cal F}_{p,q,\alpha}\left(\frac{\widehat{g}}{\|\widehat{g}\|_{W^1_2(\Omega)}}\right)\ge-\lambda_{p,q,\alpha}\|\widehat{g}\|^2_{W^1_2(\Omega)}.
\end{equation*}
Since $\lambda_{p,q,\alpha}>0$ and $\|\widehat{g}\|_{W^1_2(\Omega)}\le1$, this gives $\|\widehat{g}\|_{W^1_2(\Omega)}=1$. So we have strong convergence, and $\widehat{g}$ provides the minimum of ${\cal F}_{p,q,\alpha}$ on $\mathcal{W}_1$.\hfill$\square$}\medskip

From now on, we consider only $\alpha$ for which $\lambda_{p,q,\alpha}\ge0$ is attained. For such $\alpha$ we denote by $g$ the minimizer of ${\cal F}_{p,q,\alpha}$ on $\mathcal{W}_1$. The Euler-Lagrange equation for $g$ is
\begin{multline*}
-(p-2)\,\mathrm{div}\Bigl(|\nabla\valpha|^{p-4}\langle\nabla\valpha,\nabla g\rangle
\nabla\valpha\Bigr)-\mathrm{div}\Bigl(|\nabla\valpha|^{p-2}\nabla g \Bigr)+(p-1)\valpha^{p-2}g \nonumber\\
 =(q-1)\mu_{p,q,\alpha}^{q/p}|x|^{\alpha}\valpha^{q-2}g+\lambda\Delta g-\lambda g.
\end{multline*}
Multiplying by $g$ and integrating we get $\lambda=\lambda_{p,q,\alpha}$. 

This equation admits the separation of variables, so we can write $g(r,\theta)=\mathbf{h}(r)\psi(\theta)$, where $\psi$ is the first spherical harmonic (cf. the proof of Lemma~\ref{lm:2}) and $\mathbf{h}$ satisfies the equation
\begin{align}
&-(p-1)\left(r^{n-1}(\valpha')^{p-2}\mathbf{h}'\right)'+(n-1)r^{n-3}(\valpha')^{p-2}\mathbf{h}+(p-1)r^{n-1}\valpha^{p-2}\mathbf{h} \nonumber\\
&\qquad =(q-1)\mu_{p,q,\alpha}^{q/p}r^{n+\alpha-1}\valpha^{q-2}\mathbf{h}+\lambda_{p,q,\alpha}\bigl(r^{n-1}\mathbf{h}'\bigr)'-(n-1)\lambda_{p,q,\alpha}r^{n-3}\mathbf{h}-\lambda_{p,q,\alpha}r^{n-1}\mathbf{h},
\label{eq:14}
\end{align}
with natural boundary condition
\begin{equation}
\label{eq:15}
\lim_{r\to1-0}\left((p-1)(\valpha'(r))^{p-2}+\lambda_{p,q,\alpha}\right)\mathbf{h}'(r)=0.
\end{equation}
Since ${\cal F}_{p,q,\alpha}(|\mathbf{h}|\psi)={\cal F}_{p,q,\alpha}(\mathbf{h}\psi)$, we may assume $\mathbf{h}\ge0$. Some additional properties of the function $\mathbf{h}$ are proved in Lemma~\ref{lm:main} (see Appendix).
\medskip

In the following two Theorems we prove inequality~\eqref{eq:D2} under assumptions of Theorem~\ref{th:main}, part 1 and part 2, respectively. 

\begin{Theorem}
\label{q_loc}
Let $2 <p<n$. Then there exists $q_{loc}(n,p)>p$ such that for each $q\in(p,q_{loc})$ and for all sufficiently large $\alpha$, we have $D^2{\cal Q}_{p,q,\alpha}(\valpha;g,g)>0$ for any non-zero $g\in \mathcal{W}$. 
\end{Theorem}

\Dok{Recall that it suffices to prove ${\cal F}_{p,q,\alpha}(g)>0$ for all non-zero $g\in{\cal W}_1$. Suppose that for some $q>p$ there exists a sequence $\alpha_k\to\infty$ such that ${\cal F}_{p,q,\alpha_k}$ takes non-positive values. By Lemma~\ref{lm:achiev}, there exist minimizers $g_k(r,\theta)=\mathbf{h}_k(r)\psi(\theta)\in{\cal W}_1$ with ${\cal F}_{p,q,\alpha_k}(g_k)\le0$. Without loss of generality, assume $g_k\rightharpoondown g(r,\theta)=\mathbf{h}(r)\psi(\theta)$ in $W^1_2(\Omega)$, $g_k\to g$ in $L_2(\Omega)$, and $\mathbf{h}_k(1)\to \mathbf{h}(1)$.

If $\mathbf{h}(1)=0$, then Lemma~\ref{lm:main} (see Appendix) gives
$$
\|g_k\|^2_{W^1_2(\Omega)}=\int\limits_S\psi^2dS\int\limits_0^1\bigl((\mathbf{h}_k')^2 r^{n-1}+(n-1)\mathbf{h}_k^2r^{n-3}+\mathbf{h}_k^2r^{n-1}\bigr)dr
\le C\int\limits_0^1 r^{n-3}dr\cdot \mathbf{h}_k^2(1)\to0,
$$
contradicting $g_k\in{\cal W}_1$. Hence $\mathbf{h}(1)\ne0$, so $g(1,\theta)\not\equiv0$.

We obtain
\begin{multline*}
0\ge\liminf{\cal F}_{p,q,\alpha_k}(g_k)\stackrel{(*)}{\ge}(p-2)\int\limits_\Omega|\nabla\varphi_p|^{p-4}\langle\nabla\varphi_p,\nabla g\rangle^2dx
+\int\limits_\Omega|\nabla\varphi_p|^{p-2}|\nabla g|^2dx\\
+(p-1)\int\limits_\Omega\varphi_p^{p-2}g^2dx-(q-1)({\rm meas}\,S)^{\frac2p-1}\lambda_p^{\frac2p}\int\limits_S g^2(1,\theta)dS\\
\stackrel{(**)}{=}\left(F_p(g)-(q-1)({\rm meas}\,S)^{\frac2p-1}\lambda_p^{\frac2p}\right)\int\limits_S g^2(1,\theta)dS
\stackrel{(***)}{\ge}K(n,p,q)\int\limits_S g^2(1,\theta)dS,
\end{multline*}
where
\begin{equation}
K(n,p,q):=\lambda_p^{\frac2p}({\rm meas}\,S)^{\frac2p-1}\Bigl(\lambda_p^{-\frac{p}{p-1}}(1-(n-1)\lambda_p)-(q-1)\Bigr).
\label{eq:K_npq}    
\end{equation}
Inequality $(*)$ follows from Lemmas~\ref{prop:lim_S} and \ref{prop:lim_inf} (see Appendix), $(**)$ is definition~\eqref{eq:Fpg}, and inequality $(***)$ is proved in Lemma~\ref{lm:2}.  

By Lemma~\ref{lm:main_ineq} (see Appendix), we have $K(n,p,p)>0$. Thus there exists $q_{loc}>p$ such that $K(n,p,q)>0$ for all $q\in(p,q_{loc})$, contradicting the inequality above, and Theorem follows.
\hfill$\square$}

 \begin{Theorem}
\label{th:p_loc}
Let $n\ge4$. Then there exists $p_{loc}(n)\in(2,n)$ such that for any $p\in(2,p_{loc})$ and $q\in(p,p^*)$, for all sufficiently large $\alpha$ we have $D^2{\cal Q}_{p,q,\alpha}(\valpha;g,g)>0$ for any non-zero $g\in \mathcal{W}$.
\end{Theorem}

\Dok{From the proof of Theorem~\ref{q_loc} we see that it suffices to prove that $K(n,p,q)>0$ for mentioned values of $p$ and $q$. Since $K(n,p,q)$ is monotonically decreasing in $q$, it suffices to prove $K(n,p,p^*)>0$ for $p>2$ close to $2$. 

The latter inequality is equivalent to $\lambda_p(n-1)+\lambda_p^{\frac p{p-1}}(p^*-1)<1$. Since $\lambda_p<\frac 1n$ (see~\eqref{eq:1n}), it suffices to prove
$$
\frac{n-1}{n}+n^{-\frac p{p-1}}(p^*-1)<1
\qquad\Longleftrightarrow\qquad n^{\frac 1{p-1}}>p^*-1=
\frac{np}{n-p}-1.
$$
If $n\ge 4$ then this inequality holds for $p=2$ (see~\cite[Prop. 3.6]{GazSer}), and therefore for all $p$ in some right neighborhood of $2$. \hfill$\square$ \medskip
}

\Zam{For $n=3$, the statement of Theorem~\ref{th:p_loc} fails. Indeed, for $p=2$, the Steklov eigenvalues are known exactly~\cite{LDT}:
$$
\lambda_p|_{p=2}=1-\frac n2+\frac{I'_{\frac n2-1}(1)}{I_{\frac n2-1}(1)},
$$
where $I_\nu$ is the modified Bessel function. For $n=3$, we have $\lambda_p|_{p=2}\simeq0.313035$, and one finds $K(3,2,2^*)<0$. This implies that $D^2{\cal Q}_{2,q,\alpha}(\valpha;g,g)$ is sign-changing for $q$ close to $2^*$ and $\alpha$ large (see~\cite[Prop. 3.6 and Rem. 3.7]{GazSer}).}

\section{Appendix}

\begin{Lemma}
\label{lm:main}
Let $\mathbf{h}_{\alpha}\in W^1_2\bigl((0,1),r^{n-1}dr\bigr)$ be nonnegative solutions of~\eqref{eq:14}--\eqref{eq:15}. Then for sufficiently large $\alpha$ we have:
\begin{itemize}
\item[$(a)$] as $r\to+0$,
$$
\mathbf{h}_{\alpha}'(r)\asymp\begin{cases}
r^{-\frac{p-2}{p-1}} & \text{if }\lambda_{p,q,\alpha}=0,\\
1 & \text{if }\lambda_{p,q,\alpha}>0.
\end{cases}
$$
\item[$(b)$] $\mathbf{h}_{\alpha}'(r)>0$ for all $r\in(0,1)$.

\item[$(c)$] the following inequality holds for $r\in(0,1):$  
\begin{equation}
\label{eq:K}
\mathbf{h}_{\alpha}'(r)\le\frac{C}{r}\,\mathbf{h}_{\alpha}(1).
\end{equation}
\end{itemize}
\end{Lemma}

\Dok{We rewrite~\eqref{eq:14} as
\begin{equation}
\label{eq:18}
\Bigl(r^{n-1}\bigl((p-1)(\valpha')^{p-2}+\lambda_{p,q,\alpha}\bigr)\mathbf{h}_{\alpha}'\Bigr)'
=V_{p,q,\alpha}(r)r^{n-1}\mathbf{h}_\alpha,
\end{equation}
where
\begin{equation}
\label{eq:def_V}
V_{p,q,\alpha}(r)=\frac{(n-1)\bigl((\valpha')^{p-2}+\lambda_{p,q,\alpha}\bigr)}{r^2}+(p-1)\valpha^{p-2}+\lambda_{p,q,\alpha}
-(q-1)\mu_{p,q,\alpha}^{q/p}r^{\alpha}\valpha^{q-2}.
\end{equation}

Due to~\eqref{eq:varphi} and~\eqref{eq:v_alpha_prime} we have, as $r\to+0$,
\begin{equation*}
V_{p,q,\alpha}(r)=(n-1)\Bigl(\varphi_0^{p-2}n^{-\frac{p-2}{p-1}}r^{-\frac{p}{p-1}}\bigl(1+o_r(1)\bigr)\bigl(1+o_\alpha(1)\bigr)+\lambda_{p,q,\alpha}r^{-2}\Bigr) +C+\lambda_{p,q,\alpha}+o_r(1)    
\end{equation*}
and
\begin{multline}
\label{eq:19}
\Bigl(r^{n-1}\bigl((p-1)n^{-\frac{p-2}{p-1}}r^{\frac{p-2}{p-1}}(1+o_r(1))(1+o_\alpha(1))+\lambda_{p,q,\alpha}\bigr)\mathbf{h}_{\alpha}'\Bigr)' \\
 =(n-1)\Bigl(n^{-\frac{p-2}{p-1}}r^{-\frac{p}{p-1}}\bigl(1+o_r(1)\bigr)\bigl(1+o_\alpha(1)\bigr)+\lambda_{p,q,\alpha} r^{-2}\Bigr)r^{n-1}\mathbf{h}_\alpha +\Bigl(C+\lambda_{p,q,\alpha}+o_r(1)\Bigr) r^{n-1}\mathbf{h}_\alpha.
\end{multline}
We seek the asymptotics of the solution in the form
$$
\mathbf{h}_\alpha(r)\sim r^{\beta}.
$$

If $\lambda_{p,q,\alpha}=0$ then we have
$$
\Bigl((p-1)r^{n-1+\frac{p-2}{p-1}}\beta r^{\beta-1}\Bigr)'\sim
(n-1)r^{n-1-\frac{p}{p-1}} r^{\beta},\qquad r\to +0;
$$
$$
(p-1)\beta\Bigl(n-1+\dfrac{p-2}{p-1}+\beta-1\Bigr)=n-1;
$$
$$
\beta=\dfrac{1}{p-1}\quad\mbox{or\quad}\beta=-(n-1).
$$
The second solution does not belong to $W^1_2((0,1),r^{n-1}dr)$. Therefore, $\mathbf{h}_\alpha\asymp r^{\frac{1}{p-1}}$. Integrating~\eqref{eq:19} from $0$ to $r$, we obtain
$$
r^{n-1+\frac{p-2}{p-1}} \mathbf{h}_\alpha'(r)\asymp r^{n-1-\frac{p}{p-1}+\frac{1}{p-1}+1},\qquad \mathbf{h}_\alpha'(r)\asymp r^{-\frac{p-2}{p-1}}. 
$$

If $\lambda_{p,q,\alpha}>0$ then we have 
$$
\Bigl(r^{n-1}\beta r^{\beta-1}\Bigr)'\sim
(n-1)r^{n-3} r^{\beta},\qquad r\to +0;
$$
$$
\beta(n+\beta-2)=n-1;\qquad\beta=1\quad\mbox{or\quad}\beta=-(n-1).
$$
The second solution does not belong to $W^1_2((0,1),r^{n-1}dr)$. Therefore, $\mathbf{h}_\alpha\asymp r$. Integrating~\eqref{eq:19} from $0$ to $r$, we obtain
\begin{multline*}
r^{n-1}\bigl(Cr^{\frac{p-2}{p-1}}(1+o_r(1))(1+o_\alpha(1))+\lambda_{p,q,\alpha}\bigr)\mathbf{h}_{\alpha}'(r)\\
 =\Bigl(Cr^{-\frac{p}{p-1}}\bigl(1+o_r(1)\bigr)\bigl(1+o_\alpha(1)\bigr)+C\lambda_{p,q,\alpha} r^{-2}\Bigr)r^{n+1}+\Bigl(C+\lambda_{p,q,\alpha}+o_r(1)\Bigr)r^{n+1},
\end{multline*}
\begin{equation*}
\mathbf{h}_{\alpha}'(r)=\dfrac{Cr^{\frac{p-2}{p-1}}\bigl(1+o_r(1)\bigr)\bigl(1+o_\alpha(1)\bigr)+C\bigl(1+o_r(1)\bigr)\lambda_{p,q,\alpha}}{Cr^{\frac{p-2}{p-1}}\bigl(1+o_r(1)\bigr)\bigl(1+o_\alpha(1)\bigr)+\lambda_{p,q,\alpha}}\asymp 1,
\end{equation*}
and the statement $(a)$ holds.\medskip

Further, statement $(a)$ together with~\eqref{eq:15} implies $r^{n-1}\bigl((p-1)(\valpha')^{p-2}+\lambda_{p,q,\alpha}\bigr)\mathbf{h}_{\alpha}'\to0$ as $r\to0$ or $r\to1-0$. Therefore, equation~\eqref{eq:18} and non-negativity of $\mathbf{h}_{\alpha}$ show that $V_{p,q,\alpha}(r)$ should be sign-changing on $(0,1)$. 

Lemma~\ref{lm:r_0} implies that it changes sign exactly once, whence $r^{n-1}\bigl((p-1)(\valpha')^{p-2}+\lambda_{p,q,\alpha}\bigr)\mathbf{h}_{\alpha}'$ increases on $(0,r_0)$ and decreases on $(r_0,1)$ for some $r_0\in(0,1)$. Thus, $\mathbf{h}_{\alpha}'>0$ on $(0,1)$ and we obtain $(b)$.
\medskip

Finally, we prove~\eqref{eq:K}. Let $r\in(0,1-\frac{M}\alpha]$, where $M$ was defined in Lemma~\ref{lm:v_alpha_1}. We integrate~\eqref{eq:18} from $0$ to $r$ using the monotonicity of $\mathbf{h}_\alpha$ and~\eqref{eq:asymp_3}:
\begin{multline*}
r^{n-1}\left((p-1)(\valpha')^{p-2}+\lambda_{p,q,\alpha}\right)\mathbf{h}_{\alpha}'(r)
=\int\limits_0^r V_{p,q,\alpha}(s)s^{n-1}\mathbf{h}_\alpha(s) ds\\
\le\int\limits_0^r\left(Cs^{\frac{p-2}{p-1}-2}+C+\lambda_{p,q,\alpha}\Bigl(1+\frac{C}{s^2}\Bigr)\right)s^{n-1}ds\cdot \mathbf{h}_\alpha(r)\le\left(Cr^{\frac{p-2}{p-1}+n-2}+C\lambda_{p,q,\alpha}r^{n-2}\right)\mathbf{h}_\alpha(r),
\end{multline*}
so
$$
\mathbf{h}_{\alpha}'(r)\le\frac{Cr^{\frac{p-2}{p-1}+n-2}+C\lambda_{p,q,\alpha}r^{n-2}}{\bigl(Cr^{\frac{p-2}{p-1}}+\lambda_{p,q,\alpha}\bigr)r^{n-1}}\,\mathbf{h}_\alpha(r)
\le\frac{C}{r}\,\mathbf{h}_\alpha(r)\le\frac{C}{r}\,\mathbf{h}_\alpha(1).
$$

Now let $r\in(1-\frac{M}\alpha,1)$. We integrate~\eqref{eq:18} from $r$ to $1$ using~\eqref{eq:15} and~\eqref{eq:asymp_1}:
\begin{multline*}
r^{n-1}\bigl((p-1)(\valpha')^{p-2}+\lambda_{p,q,\alpha}\bigr)\mathbf{h}_{\alpha}'(r)
\le\int\limits_r^1(q-1)\mu_{p,q,\alpha}^{q/p}s^{\alpha}\valpha^{q-2}(s)s^{n-1}\mathbf{h}_\alpha(s) ds\\
\le C\alpha\int\limits_r^1 s^{\alpha+n-1}ds\cdot \mathbf{h}_\alpha(1)\le C(1-r^{n+\alpha})\mathbf{h}_\alpha(1)\le C\alpha(1-r)\mathbf{h}_\alpha(1).
\end{multline*}
Then \eqref{eq:asymp_1} implies that
\begin{equation*}
\mathbf{h}_{\alpha}'(r)\le\dfrac{C\alpha(1-r)}{C\alpha^{\frac{p-2}{p-1}}(1-r)^{\frac{p-2}{p-1}}+\lambda_{p,q,\alpha}} \,\mathbf{h}_\alpha(1)\le C\alpha^{\frac{1}{p-1}}(1-r)^{\frac{1}{p-1}}\mathbf{h}_\alpha(1)\le C \mathbf{h}_\alpha(1),
\end{equation*}
and the statement $(c)$ holds. \hfill$\square$
}

\begin{Lemma}
\label{lm:r_0}
For sufficiently large $\alpha$ the function $V_{p,q,\alpha}$ defined in \eqref{eq:def_V} is positive in a neighborhood of the origin and changes sign on $(0,1)$ not more that once. 
\end{Lemma}

\Dok First, we notice that for $r\in(0,1-\frac{3\ln\alpha}{2\alpha}]$ we have by Theorem~\ref{th:lim}$(e)$
$$
(q-1)\mu_{p,q,\alpha}^{q/p}r^{\alpha}\valpha^{q-2}(r)\le C\alpha\left(1-\tfrac{3\ln\alpha}{2\alpha}\right)^\alpha\le \dfrac{C}{\alpha^{1/2}}=o_\alpha(1),
$$ 
so $V_{p,q,\alpha}(r)>0$ for sufficiently large $\alpha$. 

Next, we claim that $V_{p,q,\alpha}$ decreases on $(1-\frac{3\ln\alpha}{2\alpha},1)$. Indeed,
\begin{multline*}
V_{p,q,\alpha}'(r)= -\frac{2(n-1)\bigl((\valpha')^{p-2}+\lambda_{p,q,\alpha}\bigr)}{r^3}+\frac{(n-1)(p-2)(\valpha')^{p-3}\valpha''}{r^2}+(p-1)(p-2)\valpha^{p-3}\valpha'\\
-(q-1)\alpha\mu_{p,q,\alpha}^{q/p}r^{\alpha-1}\valpha^{q-2}-(q-1)(q-2)\mu_{p,q,\alpha}^{q/p}r^{\alpha}\valpha^{q-3}\valpha'.   
\end{multline*}

We multiply this relation by $\valpha'>0$ (see Lemma~\ref{lm:v_alpha_1}) and express the term $(\valpha')^{p-2}\valpha''$ from equation~\eqref{eq:03_1}. This gives
\begin{multline*}
V_{p,q,\alpha}'(r)\valpha'(r)=-\frac{2(n-1)\lambda_{p,q,\alpha}\valpha'}{r^3} -\frac{2(n-1)(\valpha')^{p-1}}{r^3}
\\+\frac{(n-1)(p-2)}{(p-1)r^2}\Bigl(-\dfrac{(n-1)(\valpha')^{p-1}}{r}+\valpha^{p-1}-\mu_{p,q,\alpha}^{q/p} r^\alpha\valpha^{q-1}\Bigr)+(p-1)(p-2)\valpha^{p-3}(\valpha')^2\\
-(q-1)\alpha\mu_{p,q,\alpha}^{q/p}r^{\alpha-1}\valpha^{q-2}\valpha'-(q-1)(q-2)\mu_{p,q,\alpha}^{q/p}r^{\alpha}\valpha^{q-3}(\valpha')^2\\
\stackrel{(*)}{\le}
C-C\alpha r^\alpha-C\alpha^2 r^\alpha \left(\frac{\alpha}{\ln\alpha}\right)^{\frac{1}{p-1}}(1-r)^{\frac{1}{p-1}}-C\alpha r^\alpha \left(\frac{\alpha}{\ln\alpha}\right)^{\frac{2}{p-1}}(1-r)^{\frac{2}{p-1}},
\end{multline*}
where the inequality $(*)$ follows from~\eqref{eq:asymp},~\eqref{eq:asymp_1} and Theorem~\ref{th:lim}(\textit{b,c,e}). As $r\in(1-\frac{3\ln\alpha}{2\alpha},1)$, we deduce that $\alpha^2 r^\alpha\ge {C}{\alpha^{\frac{1}{2}}}$. 
Dropping the last term, we arrive at 
\begin{equation*}
 V_{p,q,\alpha}'(r)\valpha'(r)\le 
C-C \alpha r^\alpha
-\frac{C\alpha^{\frac{1}{2}}}{{(\ln\alpha)}^{\frac{1}{p-1}}}\,\alpha^{\frac{1}{p-1}}(1-r)^{\frac{1}{p-1}}.
\end{equation*}
For $r\in(1-\frac{3\ln\alpha}{2\alpha},1-\frac{\ln\alpha}{2\alpha})$ we drop the second term and estimate
\begin{equation*}
\alpha(1-r) \asymp \ln\alpha 
\quad \Longrightarrow\quad  V_{p,q,\alpha}'(r)\valpha'(r)\le
C - C\alpha^{\frac{1}{2}}<0,
\end{equation*}
whereas for $r\in(1-\frac{\ln\alpha}{2\alpha},1)$ we drop the third term and obtain
\begin{equation*}
\alpha r^\alpha\ge C\alpha^{\frac{1}{2}} \quad\Longrightarrow\quad V_{p,q,\alpha}'(r)\valpha'(r)\le
C - C\alpha^{\frac{1}{2}}<0,
\end{equation*}
and the claim follows. The proof is complete. \hfill $\square$
\medskip

\begin{Lemma}
\label{prop:lim_S}
Assume that $2<p<n$, $p<q<p^*$ and $g_k$ are minimizers of ${\cal F}_{p,q,\alpha_k}$ on~${\cal W}_1$. Then, up to a subsequence, $g_k\rightharpoondown g$ in $W^1_2(\Omega)$ as $\alpha_k\to\infty$, and 
\begin{equation}
\gathered
\int\limits_\Omega {\bf v}_{\alpha_k}^{p-2}g_k^2dx=\int\limits_\Omega \varphi_p^{p-2}g^2dx+o_k(1);\\
(\alpha_k+n)\int\limits_\Omega|x|^{\alpha_k}{\bf v}_{\alpha_k}^{q-2}g_k^2dx
=\varphi_p^{q-2}(1)\int\limits_S g^2(1,\theta)dS+o_k(1).
\endgathered
\label{eq:2-2}    
\end{equation}
\end{Lemma}

\Dok{By Theorem~\ref{th:lim}$(b)$ we have ${\bf v}_{\alpha_k}\rightrightarrows\varphi_p$. Next, $g_k$ are normalized in $W^1_2(\Omega)$, so we can extract a subsequence such that $g_k\rightharpoondown g$ in $W^1_2(\Omega)$. By the Rellich theorem, $g_k\to g$ in $L_2(\Omega)$ and $\|g_k\|_{L_2(\Omega)}\to \|g\|_{L_2(\Omega)}$.
Therefore
$$
\left|\int\limits_\Omega {\bf v}_{\alpha_k}^{p-2}g_k^2 dx-\int\limits_\Omega \varphi_p^{p-2}g^2 dx\right|\le
\sup\limits_k\sup\limits_{\Omega}{\bf v}_{\alpha_k}^{p-2}\int\limits_\Omega \bigl|g_k^2-g^2\bigr|dx+\int\limits_\Omega \bigl|{\bf v}_{\alpha_k}^{p-2}-\varphi_p^{p-2}\bigr|g^2dx=o_k(1).
$$
For the second relation in~\eqref{eq:2-2} we notice that $(\alpha_k+n)|x|^{\alpha_k}=\mbox{div}\bigl(|x|^{\alpha_k}x\bigr)$, therefore
\begin{multline}
(\alpha_k+n)\int\limits_{\Omega}|x|^{\alpha_k}{\bf v}_{\alpha_k}^{q-2}g_k^2 dx=\int\limits_{\Omega}{\bf v}_{\alpha_k}^{q-2}g_k^2 \,\mbox{div}\bigl(|x|^{\alpha_k}x\bigr)dx\\
=\int\limits_S |x|^{\alpha_k}{\bf v}_{\alpha_k}^{q-2}g_k^2 \langle x,{\bf n}\rangle d\theta-
(q-2)\int\limits_{\Omega}|x|^{\alpha_k}{\bf v}_{\alpha_k}^{q-3}g_k^2 \langle\nabla {\bf v}_{\alpha_k},x\rangle dx-
2\int\limits_{\Omega}|x|^{\alpha_k}{\bf v}_{\alpha_k}^{q-2}g_k \langle\nabla g_k,x\rangle dx\\
={\bf v}_{\alpha_k}^{q-2}(1)\int\limits_S g_k^2(1,\theta)d\theta-
(q-2)\int\limits_{\Omega}|x|^{\alpha_k+1}{\bf v}_{\alpha_k}^{q-3}{\bf v}_{\alpha_k}'g_k^2  dx-
2\int\limits_{\Omega}|x|^{\alpha_k}{\bf v}_{\alpha_k}^{q-2}g_k \langle\nabla g_k,x\rangle dx.
\label{f37}
\end{multline}
Since ${\bf v}_{\alpha_k}(1)\to \varphi_p(1)$, by the trace embedding theorem  we obtain
$$
{\bf v}_{\alpha_k}^{q-2}(1)\int\limits_S g_k^2(1,\theta)d\theta=\varphi_p^{q-2}(1)\int\limits_S g^2(1,\theta)d\theta+o_k(1).
$$
Next, we claim that the last two integrals in~\eqref{f37} are small. Indeed, applying Theorem~\ref{th:lim} and H\"older's inequality, we have for some $s\in(2,2^*)$
\begin{align*}
\left|\int\limits_{\Omega}|x|^{\alpha_k+1}{\bf v}_{\alpha_k}^{q-3}{\bf v}_{\alpha_k}'g_k^2  dx\right|\le &\ 
C\int\limits_{\Omega}|x|^{\alpha_k+1}g_k^2 dx\le C\left(\int\limits_{\Omega} |x|^{\frac{(\alpha_k+1)s}{s-2}}\right)^{\frac{s-2}{s}}\|g_k\|^2_{L_s(\Omega)}=o_k(1);\\
\left|\int\limits_{\Omega}|x|^{\alpha_k}{\bf v}_{\alpha_k}^{q-2}g_k \langle\nabla g_k,x\rangle dx\right|\le &\  C\left(\int\limits_{\Omega} |x|^{\frac{2(\alpha_k+1)s}{s-2}}\right)^{\frac{s-2}{2s}} \| g_k\|_{L_s(\Omega)}\cdot \|\nabla g_k\|_{L_2(\Omega)}=o_k(1),
\end{align*}
and the claim follows. So the relations~\eqref{eq:2-2} are proved. \hfill$\square$}\medskip

\begin{Lemma}\label{prop:lim_inf}
Assume that $2<p<n$, $p<q<p^*$ and $g_k$ are minimizers of ${\cal F}_{p,q,\alpha_k}$ on ${\cal W}_1$, $g_k\rightharpoondown g$. Then, as $\alpha_k\to\infty$,
\begin{multline*}
\liminf\left[(p-2)\int\limits_\Omega|\nabla{\bf v}_{\alpha_k}|^{p-4}\langle\nabla{\bf v}_{\alpha_k},\nabla g_k\rangle^2dx
+\int\limits_\Omega|\nabla{\bf v}_{\alpha_k}|^{p-2}|\nabla g_k|^2dx\right]
\\ \ge(p-2)\int\limits_\Omega|\nabla\varphi_p|^{p-4}\langle\nabla\varphi_p,\nabla g\rangle^2dx
+\int\limits_\Omega|\nabla\varphi_p|^{p-2}|\nabla g|^2dx.
\end{multline*}
\end{Lemma}

\Dok{
We prove the inequality for the second term, the first one is considered similarly. We have
\begin{equation}\label{f38}
\int\limits_\Omega |\nabla{\bf v}_{\alpha_k}|^{p-2}|\nabla g_k|^2 dx=
\int\limits_\Omega \bigl(|\nabla{\bf v}_{\alpha_k}|^{p-2}-|\nabla\varphi_p|^{p-2}\bigr)|\nabla g_k|^2 dx+
\int\limits_\Omega |\nabla\varphi_p|^{p-2}|\nabla g_k|^2 dx.
\end{equation}
The functional $\int\limits_\Omega |\nabla\varphi_p|^{p-2}|\nabla g_k|^2 dx$ is convex, so 
$$
\liminf \int\limits_\Omega |\nabla\varphi_p|^{p-2}|\nabla g_k|^2 dx\ge \int\limits_\Omega |\nabla\varphi_p|^{p-2}|\nabla g|^2 dx. 
$$
We claim that the first term in~\eqref{f38} tends to zero. Indeed, take an arbitrary $\varepsilon>0$. By the trace embedding theorem, $g_k\to g$ in $L_2(S)$. Moreover, since the function $g_k$ is a minimizer of ${\cal F}_{p,q,\alpha_k}$, we rewrite it as $g_k(r,\theta)=\mathbf{h}_k(r)\psi(\theta)$ with the first spherical harmonic $\psi$. 

Due to Theorem~\ref{th:lim}$(c)$ and~\eqref{eq:K}, we can choose $0<\delta<1$ such that 
\begin{multline*}
\int\limits_{\Omega\backslash\Omega_\delta}\bigl||\nabla{\bf v}_{\alpha_k}|^{p-2}-|\nabla\varphi_p|^{p-2}\bigr||\nabla g_k|^2 dx\le
C\int\limits_{\Omega\backslash\Omega_\delta}|\nabla g_k|^2 dx\\
= C\int\limits_S dS\int\limits_{\delta}^1 (\mathbf{h}_k'\,^2 r^{n-1}\psi^2+\mathbf{h}_k^2 r^{n-3}|\nabla_\theta \psi|^2) dr
\\
\le C\int\limits_S \psi^2 dS\cdot \mathbf{h}_k^2(1) \int\limits_{\delta}^1r^{n-3}dr \le C(1-\delta)\|g\|^2_{L^2(S)}<\varepsilon.
\end{multline*}
By Theorem~\ref{th:lim}$(d)$, ${\bf v}_{\alpha_k}\to\varphi_p$ in the space ${\cal C}^1(\overline\Omega_\delta)$. Therefore,  
$$
\left|\int\limits_\Omega \bigl(|\nabla{\bf v}_{\alpha_k}|^{p-2}-|\nabla\varphi_p|^{p-2}\bigr)|\nabla g_k|^2 dx\right|\le
o_k(1)\int\limits_{\Omega_\delta}|\nabla g_k|^2 dx+\varepsilon.
$$
Since $\varepsilon$ is arbitrarily small, the claim follows. This completes the proof. \hfill$\square$
}

\begin{Lemma}
\label{lm:main_ineq}
Let $n\ge 3$. Then for any $2\le p<n$ we have $K(n,p,p)>0$, where $K(n,p,q)$ was defined in~\eqref{eq:K_npq}.
\end{Lemma}

\Dok{It is sufficient to prove the inequality
$$
\lambda_p(n-1)+\lambda_p^{\frac p{p-1}}(p-1)<1.
$$
Consider the auxiliary function $\phi_p(r)=\exp\{\varkappa p^{-1}r^{p'}\}$, where $p'=\dfrac{p}{p-1}$ and $\varkappa=n^{-\frac 1{p-1}}(p-1)$. Then $\phi_p'(r)=n^{-\frac 1{p-1}}r^{\frac 1{p-1}}\exp\{\varkappa p^{-1}r^{p'}\}$, and from~\eqref{eq:lambda_p} we derive
\begin{equation*}
\lambda_p \le\frac{\|\phi_p\|^p_{W^1_p(\Omega)}}{\|\phi_p\|^p_{L_p(S)}} 
= \left(n^{-p'}\int\limits_0^1 r^{n-1+p'}e^{\varkappa r^{p'}}dr+\int\limits_0^1 r^{n-1}e^{\varkappa r^{p'}}dr\right)\cdot\dfrac{1}{e^{\varkappa}}.
\end{equation*}
We integrate the second term by parts and change the variable $t=r^{p'}$. This gives
\begin{equation}
\lambda_p\le \frac1n\left(1-\frac{e^{-\varkappa}\varkappa}{p'}\int\limits_0^1 t^{\frac{n}{p'}}e^{\varkappa t}dt\right)=:\frac1n(1-I_{p,n}).
\label{eq:1n}
\end{equation}
Therefore,
$$
\lambda_p(n-1)+\lambda_p^{\frac p{p-1}}(p-1)\le\frac{n-1}{n}(1-I_{p,n})+\frac{\varkappa}{n}(1-I_{p,n})^{\frac{p}{p-1}},
$$
and it is sufficient to show that
\begin{equation*}
\varkappa(1-I_{p,n})^{\frac{p}{p-1}}<(n-1)I_{p,n}+1
\end{equation*}
or
\begin{equation}
\label{nerav}    
(1-I_{p,n})^p<\left(\dfrac{(n-1)I_{p,n}+1}{\varkappa}\right)^{p-1}.
\end{equation}
This inequality is proved by a numerical-analytical method in Lemma~\ref{ner_Ipn}. \hfill$\square$}

\begin{Lemma}\label{ner_Ipn}
\it Let $n\ge 3$, $2\le p\le n$. Then inequality~\eqref{nerav} holds.
\end{Lemma}
\Dok{Denote $\beta=\dfrac{n}{p'}>1$ and $f_{\beta}(t)=t^{\beta}(1-t)$. Then
$$
\int\limits_0^1 t^{\beta}e^{\varkappa t}dt=\int\limits_0^1 t^{\beta+1}e^{\varkappa t}dt+\int\limits_0^1 f_\beta(t)e^{\varkappa t}dt=
\dfrac{e^{\varkappa}}{\varkappa}-\dfrac{\beta+1}{\varkappa }\int\limits_0^1 t^{\beta}e^{\varkappa t}dt+\int\limits_0^1 f_\beta(t)e^{\varkappa t}dt,
$$
so
\begin{equation}
\left(1+\dfrac{\beta+1}{\varkappa}\right)\dfrac{p'e^{\varkappa}}{\varkappa}I_{p,n}=\dfrac{e^{\varkappa}}{\varkappa}+\int\limits_0^1 f_\beta(t)e^{\varkappa t}dt.
\label{eq:int_f_beta}   
\end{equation}

First we claim that
\begin{equation}
f_{\beta}(t)\ge e^{-1}(1-t)\left((\beta+1)t-(\beta-1)\right)\quad \mbox{for}\ \  t\in\left[\frac{\beta-1}{\beta+1},1\right].
\label{eq:f_beta}    
\end{equation}
Indeed, the function $t^{\beta}$ is convex, so 
$$
t^{\beta}\ge \left(t_0^{\beta}+\beta t_0^{\beta-1}(t-t_0)\right)\Bigr|_{t_0=\frac{\beta}{\beta+1}}= \left(\dfrac{\beta}{\beta+1}\right)^{\beta}\left((\beta+1)t-(\beta-1)\right).
$$ 
Since $\bigl(\frac{\beta+1}{\beta}\bigr)^{\beta}<e$, inequality~\eqref{eq:f_beta} follows.

Then we can estimate the integral on the right-hand side of~\eqref{eq:int_f_beta} as follows:
\begin{multline*}
\int\limits_0^1 f_\beta(t)e^{\varkappa t}dt > \dfrac{1}{e}\int\limits_{\frac{\beta-1}{\beta+1}}^1 (1-t)\left((\beta+1)t-(\beta-1)\right) e^{\varkappa t}dt\\
=\dfrac{2}{e\varkappa} \int\limits_{\frac{\beta-1}{\beta+1}}^1 \left((\beta+1)t-\beta\right) e^{\varkappa t}dt=
\left.\dfrac{2}{e}\left(\dfrac{\varkappa\left((\beta+1)t-\beta\right)-(\beta+1)}{\varkappa^3}\right)e^{\varkappa t}\right|_{\frac{\beta-1}{\beta+1}}^1\\
=\dfrac{2e^{\varkappa}}{e\varkappa^2}\left( \dfrac{\varkappa-(\beta+1)}{\varkappa}+\dfrac{\varkappa+(\beta+1)}{\varkappa}e^{-\frac{2\varkappa}{\beta+1}}\right)=:\dfrac{e^{\varkappa}}{\varkappa^2}\cdot g(\tpn),
\end{multline*}
where $\tpn=\dfrac{\varkappa}{\beta+1}=\dfrac{n^{-\frac{1}{p-1}}p}{n+p'}\in(0,1)$ and 
$$
g(t)=\dfrac{2}{et}\bigl(e^{-2t}(t+1)+t-1\bigr).
$$
It is easy to see that the expression in the brackets vanishes at $t=0$, and its derivative is positive. Whence, $g(t)>0$ on $(0,1)$. From~\eqref{eq:int_f_beta} we obtain
\begin{equation*}
I_{p,n}>\dfrac{1}{p'}\dfrac{t_{p,n}}{t_{p,n}+1}\left(1+\dfrac{g(\tpn)}{\varkappa}\right)
\end{equation*}

We substitute this estimate into both parts of~\eqref{nerav}, multiply by $(1+t_{p,n})^p$ and obtain 
\begin{align*}
(1-I_{p,n})^p(1+t_{p,n})^p\le\left(1+\dfrac{t_{p,n}}{p}-\dfrac{g(\tpn)}{n+p'}\right)^p;\\
\left(\dfrac{(n-1)I_{p,n}+1}{\varkappa}\right)^{p-1}(1+t_{p,n})^p\ge (1+\epn)^{p-1}(1+\tpn),
\end{align*}
where $\epn=\dfrac{1}{\varkappa}\left(\dfrac{n-1}{p'\varkappa}t_{p,n}\cdot g(t_{p,n})-\dfrac{t_{p,n}}{p'}+1\right)$. So, the statement of the Lemma follows from the inequality
\begin{equation}
\left(1+\dfrac{t_{p,n}}{p}-\dfrac{g(\tpn)}{n+p'}\right)^p<(1+\epn)^{p-1}(1+\tpn).\label{f35}
\end{equation}

By~\eqref{eq:1n}, $1-I_{p,n}>0$, so the bracket on the left-hand side of~\eqref{f35} is positive, and we can estimate it using the inequality  $(1+y)^p\le e^{py}$:
\begin{multline*}
\left(1+\dfrac{t_{p,n}}{p}-\dfrac{g(\tpn)}{n+p'}\right)^p\le\exp\left(t_{p,n}-\frac{pg(t_{p,n})}{n+p'}\right)\\ \le
\exp\left(t_{p,n}-\frac{pn^{-\frac{1}{p-1}}}{n+p'}\,g(t_{p,n})\right)=
\exp\bigl(\tpn(1-g(\tpn))\bigr).
\end{multline*}

To estimate the right-hand side of~\eqref{f35}, we apply Bernoulli's inequality and recall that $\frac{p-1}{\varkappa}>1$:
\begin{multline*}
(1+\epn)^{p-1}>1+(p-1)\epn=1+\dfrac{p-1}{\varkappa}\left(\dfrac{n-1}{p'\varkappa}t_{p,n}\cdot g(t_{p,n})-\dfrac{t_{p,n}}{p'}+1\right)\\
\ge 1+\left(\dfrac{n^{\frac{1}{p-1}}(n-1)}{p}t_{p,n}\cdot g(t_{p,n})-t_{p,n}+1\right).
\end{multline*}

The fraction $\frac{n^{\frac{1}{p-1}}(n-1)}{p}$ decreases with $p$, so for $2\le p\le n$ we have
$$
\dfrac{n^{\frac{1}{p-1}}(n-1)}{p}\ge n^{-\frac{n-2}{n-1}}(n-1)>1
$$
(the last inequality can be proved by elementary calculus).

Substituting these estimates into~\eqref{f35} we obtain that it suffices to prove the inequality
$$
e^{t(1-g(t))}<(2-t(1-g(t))(1+t)\qquad\forall\, t\in(0,1),
$$
or
\begin{equation}
\widetilde{G}(t):= e^{G(t)}+G(t)(t+1)<2(t+1)\qquad\forall t\in(0,1),\label{f39}
\end{equation}
where 
$$
G(t)=t(1-g(t))=t-\dfrac{2}{e}\bigl(e^{-2t}(t+1)+t-1\bigr).
$$
It is easy to see that $G'(t)=\left(1-\frac{2}{e}\right)+2e^{-2t-1}(2t+1)>0$ for all $t\in(0,1)$, so the function $\widetilde{G}(t)$ is increasing. We construct the mesh $t_k=0.05k$, $k=0,1,\ldots 20$, and~\eqref{f39} follows from the chain of inequalities
$$
\widetilde{G}(t)<\widetilde{G}(t_{k})\stackrel{(*)}{<}2(t_{k-1}+1)<2(t+1),\qquad t\in[t_{k-1},t_k],\quad k=1,\ldots 20.
$$
Inequality $(*)$ is shown in the table below.\bigskip

\begin{tabular}{|l|l|l||l|l|l|}
\hline 
$k$ &$\widetilde G(t_k)\ (\pm 10^{-4})$& $2(t_{k-1}+1)$& $k$ &$\widetilde G(t_k)\ (\pm 10^{-4})$&$2(t_{k-1}+1)$\\ 
\hline\hline 
1&1.1036&2.0&11&2.4284&3.0\\ 
\hline 
2&1.2142&2.1&12&2.5791&3.1\\ 
\hline 
3&1.3310&2.2&13&2.7316&3.2\\ 
\hline 
4&1.4536&2.3&14&2.8855&3.3\\ 
\hline 
5&1.5813&2.4&15&3.0407&3.4\\ 
\hline 
6&1.7137&2.5&16&3.1970&3.5\\ 
\hline 
7&1.8501&2.6&17&3.3541&3.6\\ 
\hline 
8&1.9903&2.7&18&3.5119&3.7\\ 
\hline 
9&2.1337&2.8&19&3.6703&3.8\\ 
\hline 
10&2.2796&2.9&20&3.8291&3.9\\ 
\hline 
\end{tabular}
}


\small

\begin{thebibliography}{99}

\bibitem{GazSer} M. Gazzini, E. Serra. The Neumann problem for the H\'enon equation, trace inequalities and Steklov eigenvalues. Ann. Inst. Henri Poincar\'e, Analyse Non Lineaire \textbf{25}, 281--302 (2008).

\bibitem{Hen} M. H\'enon. Numerical experiments on the stability of spherical stellar systems. Astronomy and Astrophysics \textbf{24}, 229--238 (1973).

\bibitem{KapLad} L.V.~Kapitanskii, O.A.~Ladyzhenskaya.
On the Coleman principle for finding stationary points of invariant functionals. J. Soviet Math., \textbf{27}(6), 2606--2616 (1984). Translated from Zap. Nauchn. Sem. LOMI, \textbf{127}, 84--102 (1983).

\bibitem{NazKol} S.B. Kolonitskii, A.I. Nazarov. Multiplicity of solutions to the Dirichlet problem for generalized H\'enon equation. J. Math. Sci. \textbf{144}(6), 4624--4644 (2007). Transl. from Problemy Mat. Analiza \textbf{35}, 91--110 (2007).

\bibitem{LDT} E. Lami Dozo, O. Torn\'e. Symmetry and symmetry breaking for minimizers in the trace inequality. Comm. Contemp. Math. \textbf{7}(6), 727--746 (2005).

\bibitem{MR} S. Mart\'inez, J.D. Rossi. Isolation and simplicity for the first eigenvalue of the $p$-Laplacian with a nonlinear boundary condition. Abstr. and Appl. Anal. \textbf{7}(5), 287--293 (2002). 

\bibitem{Naz} A.I. Nazarov. On the symmetry of extremals in the weight embedding theorem. J. Math. Sci. \textbf{107}(3), 3841--3859 (2001). Transl. from Problemy Mat. Analiza \textbf{23}, 50--75 (2001).

\bibitem{Ni} W.-M. Ni. A nonlinear Dirichlet problem on a unit ball and its applications. Indiana Univ. Math. J. \textbf{31}(6), 801--807 (1982).

\bibitem{Pal} R.S.~Palais. The principle of symmetric criticality. Comm. Math. Phys. \textbf{69}(1), 19--30 (1979).

\bibitem{Shch18} A.P. Shcheglova. The Neumann problem for the generalized H\'enon equation. J. Math. Sci. \textbf{235}(3), 360--373 (2018). Transl. from Problemy Mat. Analiza \textbf{95}, 3--19 (2018).

\bibitem{SSW} D. Smets, M. Willem, J. Su. Non-radial ground states for the H\'enon equation. Commun. Contemp. Math. \textbf{4}(3), 467--480 (2002).

\bibitem{Str} W. A. Strauss. Existence of solitary waves in higher dimensions. Commun. Math. Physics \textbf{55}, 149--162 (1977).

\end{thebibliography}
\end{document}